\documentclass[10pt]{article}
\usepackage[cp866]{inputenc}
\begin{document}

\centerline {{\Large\bf Skew-symmetric differential forms. Invariants.}}
\centerline {{\Large\bf Realization of invariant structures.}}
\centerline {\it L.~I. Petrova}

\renewcommand{\abstractname}{Abstract}
\begin{abstract}

Skew-symmetric differential forms play an unique role in mathematics
and mathematical physics. This relates to the fact that closed exterior
skew-symmetric differential forms are invariants. 
The concept of ``Exterior differential forms" 
was introduced by E.Cartan for a notation of integrand
expressions, which can create the integral invariants.
(The existence of integral invariants was recognized by A.
Poincare while studying the general equations of dynamics.)

All invariant mathematical formalisms are based on invariant properties 
of closed exterior forms. The invariant properties of closed exterior 
forms explicitly or implicitly manifest themselves essentially in all
formalisms of field theory, such as the Hamilton formalism, tensor
approaches, group methods, quantum mechanics equations, the Yang-Mills
theory and others. They lie at the basis of field theory.

However, in this case the question of how the closed exterior
forms are obtained arises. In present work it is shown that closed
exterior forms, which possess the invariant properties, are
obtained from skew-symmetric differential forms, which, as
contrasted to exterior forms, are defined on nonintegrable
manifolds. The process of generating closed exterior forms
describes the mechanism of realization of invariants and invariant
structures.

\end{abstract}

\section{Closed exterior skew-symmetric differential forms: Invariants.
Invariant structures.}

Distinguishing properties of the mathematical apparatus of
exterior differential forms were formulated by Cartan [1]: ``\dots
I wanted to build the theory, which contains concepts and
operations being {\it independent of any change of variables both
dependent and independent}; to do so it is necessary to change
{\it partial derivatives by differentials} that have {\it
interior meaning}.''

\subsection{Some foundations of closed exterior differential
forms}

The exterior differential form of degree $p$ ($p$-form) on integrable
manifold  can be written as [2,3]
$$
\theta^p=\sum_{i_1\dots i_p}a_{i_1\dots i_p}dx^{i_1}\wedge
dx^{i_2}\wedge\dots \wedge dx^{i_p}\quad 0\leq p\leq n\eqno(1)
$$
Here $a_{i_1\dots i_p}$ are functions of variables $x^{i_1}$,
$x^{i_2}$, \dots, $x^{i_n}$, $n$ is the dimension of space,
$\wedge$ is the operator of exterior multiplication, $dx^i$,
$dx^{i}\wedge dx^{j}$, $dx^{i}\wedge dx^{j}\wedge dx^{k}$, \dots\
is the local basis which satisfies the condition of exterior
multiplication:
$$
\begin{array}{l}
dx^{i}\wedge dx^{i}=0\\
dx^{i}\wedge dx^{j}=-dx^{j}\wedge dx^{i}\quad i\ne j
\end{array}
$$
[In further presentation the symbol of summing $\sum$ and the symbol
of exterior multiplication $\wedge$ will be omitted. Summation
over repeated indices is implied.]

An exterior differential form is called a closed one
if its differential is equal to zero:
$$
d\theta^p=0\eqno(2)
$$
From condition (2) one can see that the closed form is a conservative
quantity. This means that this can correspond to the conservation law,
namely, to some conservative physical quantity.

The differential of the form is a closed form. That is
$$
dd\omega=0
$$
where $\omega$ is an arbitrary exterior form.

The form which is the differential of some other form:
$$
\theta^p=d\omega\eqno(3)
$$
is called the exact form. The exact forms prove to be closed
automatically
$$
d\theta^p=dd\omega=0\eqno(4)
$$

Here it is necessary to pay attention to the following points. In the
above presented formulas it was implicitly assumed that the differential
operator $d$ is a total one (that is, the
operator $d$ acts everywhere in the vicinity of the point
considered). However, the differential may be internal. Such a
differential acts on some structure with the dimension being less
than that of the initial manifold.

If the exterior form is closed only on structure, the closure
condition is written as
$$
d_\pi\theta^p=0\eqno(5)
$$
In this case the structure $\pi$ obeys the condition
$$
d_\pi{}^*\theta^p=0\eqno(6)
$$
where ${}^*\theta^p$ is the dual form.

Such an exterior form is called the closed inexact form.

The structure, on which the exterior
differential form may become a closed (inexact) form, is
a pseudostructure with respect to its metric properties.

From conditions (5) and (6) one can see that the form
closed on pseudostructure is a conservative object, namely, this
quantity conserves on pseudostructure. This can also correspond to
some conservation law, i.e. to conservative object.

\subsection*{Pseudostructures}

As one can see from condition (6), the structure, on which a closed 
(inexact) form is defined, is described by dual form. The dual form is 
a closed metric form of this structure.

To understand the properties of such structure, one can use the 
correspondence between the exterior differential form and skew-symmetric 
tensor. It is known that the skew-symmetric tensors
correspond to closed exterior differential forms, and the pseudotensors 
correspond to relevant dual forms. This points to the fact that 
the structures, on which closed inexact forms are defined, are 
pseudostructures.

The characteristics, integral surfaces, surfaces
of potential (of simple layer, double layer), sections of cotangent
bundles (Yang-Mills fields), cotangent
manifold, eikonals, cohomologies by de Rham, singular cohomologies,
the pseudo-Riemann and pseudo-Euclidean spaces and others can be
regarded as examples of pseudostructures and pseudospaces, on which
closed inexact forms are defined.

It should be emphasized that the pseudostructure and corresponding 
closed inexact form make up a differential - geometrical structure. 
As it will be shown below, such a differential - geometrical
structure proves to be an invariant structure.

\subsection*{Differentials}

The exact form is, by definition, a differential (see condition (3)).
In this case the differential is total. The
closed inexact form is a differential too. And in this case the
differential is an interior one defined on pseudostructure. Thus, any
closed form is a differential. The exact form is a total
differential. The closed inexact form is an interior (on
pseudostructure) differential, that is
$$
\theta^p_\pi=d_\pi\omega\eqno(7)
$$
At this point it is worth noting that the total differential of the
form closed on pseudostructure is nonzero, that is
$$
dd_\pi\omega\ne0\eqno(8)
$$

\subsection{Invariants. Invariant structures. }

Since the closed form is a differential (a total one if the form is 
exact, or an interior one on the pseudostructure if the form is inexact),
it is obvious that the closed form proves to be invariant under all
transformations that conserve the differential. The unitary 
transformations (0-form), the tangent and canonical transformations 
(1-form), the gradient and gauge transformations (2-form) and so on are 
examples of such transformations.

{\it These are gauge transformations for spinor, scalar, vector, tensor
fields}. It can be pointed out
that just such transformations are used in field theory.

As mentioned above, from the closure conditions it follows
that the closed form is a conservative quantity. As the result,
the closed form is a conservative invariant quantity. This property
of closed forms plays an essential role in describing the
conservation laws and lies at the basis of field theory. The
covariance of dual form is directly connected with the invariance
of exterior closed inexact form.

\subsection*{Invariant structures}

The closed inexact exterior forms are of most significance in 
mathematical formalisms and mathematical physics. This is due to the 
fact that the closed exterior form and relevant dual form describe
the differential-geometrical structure, which is invariant one.

From the definition of closed inexact exterior form one can see
that to this form there correspond two conditions:

(1) condition (5) is a closure condition of exterior form
itself, and

(2) condition (6) is that of dual form.

Conditions (5) and (6) can be regarded as equations for a binary
object that combines the pseudostructure (dual form) and the
conservative quantity (the exterior differential form) defined on
this pseudostructure. Such a binary object is a differential -
geometrical structure.  (The well-known G-Structure is an example
of such differential-geometrical structure.)

As it has been already pointed out, the closed inexact exterior form is
a differential (an interior one on pseudostructure), and hence it remains 
invariant under all transformations that conserve the differential.
Therefore, the relevant differential-geometrical structure also remains
invariant under all transformations that conserve
differential. For the sake of convenience in subsequent
presentation such differential - geometrical structures
will be called the Inv. Structures.

To an unique role of such invariant structures in mathematics it points
the fact that the transformations conserving the differential
(unitary, tangent, canonical, gradient and gauge ones) lie at the basis
of many branches of mathematics, mathematical physics and field theory.
The differential-geometrical structures made up of characteristics and
integral curves of differential equations and relevant
conditions on those are examples of Inv. Structures.

As it will be shown in Section 3 of present paper, the Inv.
Structures are of unique importance in mathematical physics and
field theory. The physical structures, of which physical fields
are made up, are such invariant structures.

It should be emphasized ones more that the Inv. Structure is a
differential-geometrical structure. That is not a spatial
structure. The spatial structure is described by {\it exact }
exterior form, whereas the invariant structure is described by
{\it inexact } exterior form.

\subsection{Invariance as the result of conjugacy of elements of
exterior or dual forms}

The closure of exterior differential forms, and hence their
invariance, results from the conjugacy of elements of exterior or
dual forms.

From the definition of the exterior differential form
one can see that exterior differential forms have complex structure.
The specific features of the exterior form structure are the homogeneity with
respect to the basis, skew-symmetry, the integration of terms each of 
which made up by two objects of different nature (the algebraic nature
of the form coefficients, and the geometric nature of the base components).
Besides,  the exterior form depends on the space dimension and on the
manifold topology. The closure property of exterior form means that any 
objects, namely, elements of exterior form, components of elements, 
elements of the form differential, exterior and dual forms and others, 
turn out to be conjugated. The variety of objects of conjugacy leads to
the fact that closed forms can describe a great number of different
invariant structures.

{\footnotesize [Let us consider some types of conjugacy.

One of the types of conjugacy is that for the form coefficients.

Let us consider the exterior differential form of first degree
$\omega=a_i dx^i$. In this case the differential will be expressed
as $d\omega=K_{ij}dx^i dx^j$, where
$K_{ij}=(\partial a_j/\partial x^i-\partial a_i/\partial x^j)$ are
the components of the form commutator.

It is evident that the differential may vanish if the components
of commutator vanish. One can see that
the components of commutator $K_{ij}$ may vanish if derivatives of
the form coefficients vanish. This is a trivial
case. Besides, the components $K_{ij}$ may vanish if the
coefficients $a_i$ are derivatives of some function $f(x^i)$,
that is, $a_i=\partial f/\partial x^i$. In this case,
the components of commutator are equal to the difference of mixed
derivatives
$$
K_{ij}=\left(\frac{\partial^2 f}{\partial x^j\partial
x^i}-\frac{\partial^2 f}{\partial x^i\partial x^j}\right)
$$
and therefore they vanish. One can see that the form
coefficients $a_i$, that satisfy these conditions,  are
conjugated quantities (the operators of mixed differentiation turn out
to be commutative).

Let us consider the case when the exterior form is written as
$$
\theta=\frac{\partial f}{\partial x}dx+\frac{\partial f}{\partial
y}dy
$$
where $f$ is the function of two variables $(x,y)$. It is evident that 
this form is closed because it is
equal to the differential $df$. And for the dual form
$$
{}^*\theta=-\frac{\partial f}{\partial y}dx+\frac{\partial
f}{\partial x}dy
$$
be also closed,  it is necessary that its commutator be equal to zero
$$
\frac{\partial^2 f}{\partial x^2}+\frac{\partial^2 f}{\partial
y^2}\equiv \Delta f=0
$$
where $\Delta$ is the Laplace operator. As a result the function $f$ has
to be a harmonic one.

Assume the exterior differential form of first degree has
the form $\theta=udx+vdy$, where $u$ and $v$ are the functions of
two variables $(x,y)$. In this case, the closure condition of the
form, that is, the condition under which the form commutator vanishes,
takes the form
$$
K=\left(\frac{\partial v}{\partial x}-\frac{\partial u}{\partial
y}\right)=0
$$
One can see that this is one of the Cauchy-Riemann conditions for
complex functions. The closure condition of the relevant dual
form ${}^*\theta=-vdx+udy$ is the second Cauchy-Riemann condition.
\{Here one can see the connection between exterior differential form and the
functions of complex variables. If we consider the function $w=u+iv$ of
complex variables $z=x+iy$ and $\overline{z}=x-iy$ that obeys the
Cauchy-Riemann conditions, then the closed exterior and dual
forms will correspond to this function. (The Cauchy-Riemann conditions 
are conditions under which the function of complex variables does not
depend on the conjugated coordinate $\overline{z}$). And
the closed exterior differential form, whose coefficients $u$ and $v$ 
are conjugated harmonic functions, corresponds to the
harmonic function of complex variables\}.

It can exist the conjugacy that makes the interior differential
on pseudostructure equal to zero, $d_\pi\theta=0$. Assume the interior 
differential is the first degree form (the form itself is a form of zero 
degree), and it can be presented as $d_\pi\theta=p_x dx+p_y dy=0$, where 
$p$ is the form of zero degree (a certain function). In this case 
the closure condition of the form is
$$
\frac{dx}{dy}=-\frac{p_y}{p_x}\eqno(9)
$$
This is a conjugacy of the basis and derivatives of the form
coefficients. One can see that this formula is one of the
formulas of canonical relations. The second formula of
canonical relations follows from the condition that the dual form
differential vanishes. This type of conjugacy  is
connected with canonical transformation. For the differential of the
first degree form (in this case the differential is a form of
second degree) the corresponding transformation has to be a gradient
transformation.
At this point it should be remarked that relation (9) is
the condition of existence of implicit function. That is, the
closed (inexact) form of zero degree is an implicit function.]}

\subsection{Identical relations of exterior differential forms, 
description of conjugacy and invariance}

Since the conjugacy is a certain connection between two operators or
mathematical objects, it is evident that the relations can  be used to
express conjugacy mathematically. Identical relations of exterior
differential forms disclose also the properties of Inv. Structure.

At this point it should be emphasized the following. The relation is a
comparison, i.e. a correlation of two objects. The relation may be
identical or nonidentical. The basis of mathematical apparatus of
exterior differential forms is made up of identical relations.
(Below nonidentical relations
will be presented, and it will be shown that identical relations for
exterior differential forms are obtained from nonidentical relations.
Also it will be shown that transitions from nonidentical relations to
identical ones describe the realization of invariant structures.)

The identical relations of exterior differential forms reflect the
closure conditions of differential forms, namely, vanishing the form
differential (see formulas (2), (5) and (6)) and the condition that
the closed differential form is a differential (see formulas (3) and 
(7)). All these conditions are the expression of conjugacy and 
invariance.

One can distinguish several types of identical relations.

1. {\it Relations in differential forms}.

They correspond to formulas (3) and (7). The examples of
such identical relations are

(a) the Poincare invariant $ds\,=\,-H\,dt\,+\,p_j\,dq_j$,

(b) the second principle of thermodynamics $dS\,=\,(dE+p\,dV)/T$,

(c) the vital force theorem in theoretical mechanics: $dT=X_idx^i$
where $X_i$ are the components of potential force, and $T=mV^2/2$ is the
vital force,

(d) the conditions on characteristics in the theory of differential
equations.

The requirement that the function is an antiderivative (the integrand is a
differential of a certain function) can be written in terms of such an
identical relation.

The existence of harmonic function is written by means
of identical relation: the harmonic function is a closed
form, that is, a differential (a differential on the Riemann surface).

In general form such an identical relation can be written as
$$
d\phi=\theta^p\eqno(10)
$$
In this relation the form in the right-hand side has to be a {\it closed }
one.

As it will be shown below,
the identical relations are satisfied only on pseudostructures. That is,
the identical relation can be written as
$$
d _{\pi}\phi=\theta _{\pi}^p\eqno(11)
$$

Identical relations (10) and (11) are the proof that the closed exterior
form is a differential, and hence, this form is an invariant with 
respect to all transformations that conserve the differential.

Identical relations occur in various branches of mathematics
and mathematical physics. Identical relations can be of another type,
namely, integral, tensor and others. And all identical relations are 
an analog to the identical relation in differential forms.

All identical relations correspond to invariant structures.

It would be noted some another types of identical relations.

2. {\it Integral identical relations}.

At the beginning of the paper it was pointed out that the exterior
differential forms were introduced as integrand expressions possessing
the following property: they can have integral invariants. This fact
(the availability of integral invariant) is mathematically expressed
as a certain identical relation.

The formulas by Newton, Leibnitz and Green, the integral relations by
Stokes and Gauss-Ostrogradskii are examples of integral identical 
relations.

3. {\it Tensor identical relations}.

From the relations that connect exterior forms of consequent degrees
one can obtain the vector and tensor identical relations that connect
the operators of gradient, curl, divergence and so on.

From the closure conditions of exterior and dual forms one can obtain
the identical relations such as the gauge relations in electromagnetic
field theory, the tensor relations between connectednesses and their
derivatives in gravitation (the symmetry of connectednesses 
with respect to lower indices,
the Bianchi identity, the conditions imposed on the
Christoffel symbols) and so on.

4. {\it Identical relations between derivatives}.

The identical relations between derivatives correspond to the closure
conditions of exterior and dual forms. The examples of such relations
are the above presented Cauchi-Riemann conditions in the theory of
complex variables, the transversality condition in the calculus of
variations, the canonical relations in the Hamilton formalism, the
thermodynamic relations between derivatives of thermodynamic functions, 
the condition that the derivative of implicit function is
subject to, the eikonal relations and so on.

\bigskip
The importance of identical relations
is manifested by the fact that practically in all branches of physics,
mechanics, thermodynamics one faces such identical relations.

The functional significance of identical relations for exterior
differential forms lies in the fact that they can describe the conjugacy
of objects that have different mathematical meaning and different 
physical nature. This enables one to see internal connections between
various branches of mathematics and physics. Due to
these possibilities the exterior differential forms, and correspondingly,
the Inv. Structures, have wide application in various
branches of mathematics and mathematical physics.

\bigskip
Identical relations possess the duality that discloses the significance
of invariant structures. The availability of differential in the 
left-hand side points to the availability of potential or state function, 
and the availability of closed inexact form points to that there is an 
invariant structure. Below it will be shown that such a relation
has a deep physical sense.

\section{Realization of invariant structures}

The mechanism of realization of invariant structures is described by
skew-symmetric differential forms, which, in contrast to
exterior forms, are defined on deforming nonintegrable manifolds 
(see Appendix of work [4]).
Such skew-symmetric differential forms possess the evolutionary 
properties. The evolutionary forms possess a peculiarity, namely,
the closed inexact exterior forms are obtained from them.
This elucidates the process of realization of invariant structures.

\subsection{Some properties of evolutionary forms}

The evolutionary skew-symmetric differential forms are obtained
from differential equations that describe any processes.

Examples of nonintegrable manifolds, on which the evolutionary
skew-symmetric differential forms are defined, are the tangent
manifolds of differential equations, the Lagrangian manifolds,
the manifolds constructed of trajectories of material medium particles
and so on. These are manifolds with unclosed metric forms. The metric 
form differential, and correspondingly its commutator, are nonzero.
(The commutators of metric forms of such manifolds describe the manifold
deformation: torsion, curvature and others).

The specific feature of evolutionary forms, i.e skew-symmetric
forms defined on deforming manifolds, is the fact that
evolutionary forms are unclosed ones. Since the basis of
evolutionary form changes, the evolutionary form differential
includes the nonvanishing differential of manifold metric form due
to differentiating the basis. Therefore, the evolutionary form
differential cannot be equal to zero. Hence, the evolutionary
form, in contrast to the case of exterior form, cannot be closed.
This leads to that in the mathematical apparatus of evolutionary
forms there arise new nonconventional elements like nonidentical
relations and degenerate transformations that allow to describe
the generation of closed inexact exterior forms and the
realization of invariant structures.

The nonidentical relations of evolutionary forms can be written as
$$
d\phi=\eta^p\eqno(12)
$$
Here $\eta^p$ is the $p$-degree evolutionary form being
unclosed, $\phi$ is some form of degree $(p-1)$, and
the differential $d\phi$ is a closed form of degree $p$.

The form differential, i.e. a closed form being an invariant object,
appears in the left-hand side of this relation. In the right-hand
side it is appeared the unclosed form, which is not an invariant
object. Such a relation cannot be identical one.

One can see the difference of relations for exterior forms and 
evolutionary ones. In the right-hand side of identical relation (see 
relation (10)) it is appeared the closed form, whereas the form in the 
right-hand side of nonidentical relation (12) is an unclosed one.

Nonidentical relations are obtained while describing any processes.
A relation of such type is obtained while, for example, analyzing the
integrability of the partial differential equation. The equation is
integrable if it can be reduced to the form $d\phi=dU$. However, it
appears that, if the equation is not subject to an additional
condition (the integrability condition), it is reduced to the
form (12), where $\eta^p$ is an unclosed form and it cannot be
written as a differential.

Nonidentical relations of evolutionary forms are evolutionary relations
because they include the evolutionary form.
Such nonidentical evolutionary relations appear to be selfvarying
ones.  The variation of any object of the relation in some
process leads to variation of another object and, in turn, the variation
of the latter leads to variation of the former. Since one of the
objects is a noninvariant (i.e. unmeasurable)
quantity, the other cannot be compared with the first one, and hence,
the process of mutual variation cannot be completed.

The nonidentity of evolutionary relation is connected with
the nonclosure of evolutionary form, that is, it is connected with
the fact that the evolutionary form commutator is nonzero.
The evolutionary form commutator includes two
terms. The first term specifies the mutual variations of evolutionary form
coefficients, and the second term (the metric form commutator) specifies
the manifold deformation. These terms have a different nature and cannot
make the commutator to be vanishing. In the process of selfvariation of
nonidentical evolutionary relation  the exchange between the
terms of evolutionary relation proceeds and this is realized according to
the evolutionary relation. The evolutionary form commutator describes
the quantity that is a moving force of evolutionary process and
leads to realization of differential-geometrical structures.

The process of the evolutionary relation selfvariation plays a governing
role in description of evolutionary processes.

The significance of the evolutionary relation selfvariation consists in
the fact that in such a process it can be realized conditions under
which the closed inexact form is obtained
from the evolutionary form and
the identical relation is obtained from the nonidentical relation.
These are conditions of degenerate transformation.
Since the evolutionary form differential is nonzero, whereas the closed
exterior form differential is zero, the transition from the evolutionary
form to closed exterior form is allowed only under {\it degenerate
transformation}. The conditions of vanishing the dual form differential
are conditions of degenerate transformation.

These are such conditions that can be realized under selfvariation of 
the nonidentical evolutionary relation.

\subsection{Realization of closed inexact exterior form. Derivation of
invariant structures}

To obtain the differential-geometrical structure, it is necessary
to obtain the closed inexact exterior form, i.e. the form closed on
pseudostructure.

To the pseudostructure it is assigned the closed dual form
(whose differential vanishes). For this reason the transition
from the evolutionary form to closed inexact exterior form proceeds
only when the conditions of vanishing the dual form differential are
realized, in other words, when the metric form differential or
commutator becomes equal to zero.

The conditions of degenerate transformation  are  conditions of
vanishing the dual form differential.
That is, it is the condition of realization of pseudostructure.
And this leads to realization of closed inexact exterior form.

As it has been already mentioned, the evolutionary differential form
$\eta^p$ involved into nonidentical relation (12) is an unclosed one.
The commutator, and hence the differential, of this form is nonzero.
That is,
$$d\eta^p\ne 0 \eqno(13)$$
If the conditions of degenerate transformation are realized, then from
the unclosed evolutionary form one can obtain the differential form 
closed on pseudostructure. The differential of this form equals zero. 
That is, it is realized the transition $d\eta^p\ne 0 \to $ (degenerate 
transformation) $\to d_\pi{}^*\eta^p=0$, $d_\pi \eta^p=0$.

The relations obtained
$$d_\pi \eta^p=0,  d_\pi{}^*\eta^p=0 \eqno(14)$$
are closure conditions for exterior inexact form, and this points to
realization of exterior form closed on pseudostructure, that is, this
points to origination of the differential-geometrical invariant 
structure.

Vanishing on pseudostructure the exterior form differential (that
is, vanishing on pseudostructure the interior differential of the
evolutionary form) points to that the exterior inexact  form is a
conservative quantity in the direction of pseudostructure.
However, in the direction normal to pseudostructure this quantity
exhibits the discontinuity. The value of such discontinuity is defined 
by the value of the evolutionary form commutator being nonzero. This
argues to discreteness of the differential-geometrical structures.

Thus, while selfvariation of the evolutionary nonidentical
relation the dual form commutator can vanish. This means that
it is made up the pseudostructure on which the differential form turns
out to be closed. The emergence of the form being closed on
pseudostructure points out to origination of invariant structures.

On the pseudostructure $\pi$ from evolutionary relation (12) it follows
the relation
$$
d_\pi\psi=\omega_\pi^p\eqno(15)
$$
which proves to be an identical relation. Indeed, since the form
$\omega_\pi^p$ is a closed one, on the pseudostructure this form turns
out to be the differential of some differential form. In other words,
this form can be written as $\omega_\pi^p=d_\pi\theta$. Relation (15)
is now written as
$$
d_\pi\psi=d_\pi\theta
$$
There are differentials in the left-hand and right-hand sides of
this relation. This means that the relation is an identical one.

From evolutionary nonidentical relation (12) it is obtained the
identical on pseudostructure relation. In this case the evolutionary
relation itself  remains to be nonidentical one. (At this point it
should be emphasized that differential, which equals zero, is an
interior one. The evolutionary form commutator becomes zero only on
pseudostructure. The total evolutionary form commutator is nonzero. That
is, under degenerate transformation the evolutionary form differential
vanishes only {\it on pseudostructure}. The total differential of
evolutionary form is nonzero. The evolutionary form remains to be
unclosed.)

It can be shown that all identical relations of the exterior
differential form theory are obtained from nonidentical relations (that
contain evolutionary forms) by applying degenerate transformations.

{\footnotesize [The conditions of degenerate transformation that lead to
origination of invariant structures can be
connected with any symmetries.  While describing material system (see,
Section 3), the symmetries can be conditioned, for example, by degrees 
of freedom of material system. Since the conditions of degenerate
transformation are those of vanishing the interior differential of
metric form, that is, vanishing the interior (rather then total)
metric form commutator, the conditions of degenerate
transformation can be caused by symmetries of coefficients of the
metric form commutator (for example, it can be the symmetric
connectedness).

Mathematically  the requirement that some functional expressions become 
equal to zero is assigned to the conditions of degenerate transformation. 
Such functional expressions are Jacobians, determinants,
the Poisson brackets, residues, and others.

The degenerate transformation is realized as the transition between 
nonequivalent frames of reference: the transition from the
noninertial frame of reference to the locally inertial one.
Evolutionary relation (12) and condition (13) are connected with the
frame of reference being related to nonintegrable noninertial manifold, 
whereas condition (14) and
identical relations (15) may be connected with only the locally inertial
frame of reference being related to pseudostructure.
For example, while studying the integrability of differential equations
under degenerate transformation it occurs the transition from the 
tangent nonintegrable  manifold to cotangent integrable  manifold. Here 
it can be underlined the connection between the degenerate 
transformation and nondegenerate one. The origination of the
differential-geometrical structures (Inv. Structures) is connected with
degenerate transformation that executes the transition from tangent
space to cotangent one. And nondegenerate transformation executes
the transition in cotangent space from any differential-geometrical
structure to another.]}

Thus, the mathematical apparatus of evolutionary differential forms
can describe the process of generation of closed inexact exterior
differential forms, and this discloses the process of origination of
invariant structures.

The process of generation of closed inexact exterior
differential forms and the origination of invariant structures
are processes of conjecting the operators. To the closed exterior form
there correspond conjugated operators, whereas to the evolutionary form
there correspond nonconjugated operators. The transition from
evolutionary form to closed exterior form and the origination of
differential-geometrical structures is a transition from nonconjugated
operators to conjugated ones. This is expressed mathematically
as the transition from nonzero differential (the evolutionary form
differential is nonzero) to the differential that equals zero (the 
closed exterior form differential equals zero).

It can be seen that the process of conjugating the objects and
obtaining the differential-geometrical structures is a mutual exchange
between the quantities of different nature (for example, between
the algebraic and geometric quantities, between the physical and
spatial quantities) and vanishing some functional expressions
(Jacobians, determinants and so on).

\subsection*{Characteristics of Inv. Structure}

Since the closed exterior differential form, which corresponds to the
Inv. Structure emerged, was obtained from evolutionary form that
enters to the nonidentical relation, it is evident that the Inv. Structure
characteristics must be connected with those of the evolutionary form and 
of the manifold on which this form is defined, as well as the conditions 
of degenerate transformation and the values of commutators of the
evolutionary form and the manifold metric form.

The conditions of degenerate transformation, as it was said before,
determine the pseudostructures. The first term of the evolutionary form
commutator determines the value of discrete change (the quantum),
which the quantity conserved on the pseudostructure undergoes under
transition from one pseudostructure to another. The second term of the
evolutionary form commutator specifies the characteristics that fixes 
the character of initial manifold deformation, which took place before
the Inv. Structure had been arisen. (Spin is such an example).

The discrete (quantum) change of a quantity proceeds in the direction
that is normal (more exactly, transverse) to the pseudostructure. Jumps
of the derivatives normal to potential surfaces are examples of such
changes.

The connection of Inv. Structure with the skew-symmetric
differential forms allows to introduce the classification of
Inv. Structures in its dependence on parameters that specify the
skew-symmetric differential forms and enter into nonidentical and
identical relation of skew-symmetric differential forms. To
determine these parameters one has to consider the problem of
integration of nonidentical evolutionary relation.

Under degenerate transformation from the nonidentical evolutionary
relation one obtains the relation being identical on pseudostructure.
Since the right-hand side of such a relation can be expressed in terms
of differential (as well as the left-hand side), one obtains the relation
that can be integrated, and as the result one obtains the relation with
differential forms of less by one degree.

The relation obtained after integration proves to be nonidentical
as well.

The resulting nonidentical relation of degree $(p-1)$ (relation that
includes the forms of degree $(p-1)$) can be integrated once again
if the corresponding degenerate transformation has been realized and
the identical relation has been formatted.

By sequential integrating the evolutionary relation of degree $p$ (in
the case of realization of corresponding degenerate transformations
and formatting the identical relation), one can get closed (on the
pseudostructure) exterior forms of degree $k$, where $k$ ranges
from $p$ to $0$.

In this case one can see that after such integration the closed (on
pseudostructure) exterior forms, which depend on two parameters, are
obtained. These parameters are the degree of evolutionary form $p$
(in the evolutionary relation) and the degree of created closed
forms $k$.

In addition to these parameters, another parameter appears, namely, the
dimension of space. If the evolutionary relation generates the closed
forms of degrees $k=p$, $k=p-1$, \dots, $k=0$, to them there correspond
the pseudostructures of dimensions $(n+1-k)$, where $n$ is the space
dimension.

\bigskip
The invariant structures are of unique significance in mathematical
physics and field theory. The physical structures that made up physical
fields are such Inv. Structures.

As it will be shown below, the mechanism of realization of 
Inv. Structures, which correspond to physical fields, describes the 
mechanism of generation of physical structures. This discloses 
the physical meaning of Inv. Structures.

\section{Physical meaning of invariant structures. Mechanism of 
generation of physical structures.}

As it has been already pointed out, the invariant structures are 
realized while analyzing the integrability of differential equations. 
Their role in the theory of differential equations relates to the fact 
that they correspond to generalized solutions which describe measurable 
physical quantities. In this case the integral surfaces with conservative 
quantities (like the characteristics, the characteristic surfaces, 
potential surfaces and so on) are invariant structures.
The examples of such studying the integrability of differential equations 
using the skew-symmetric differential forms are presented in paper [5].

\bigskip
The unique results are obtained in studying the differential equations 
that describe the conservation laws for material media. The Inv. 
Structures that correspond to
physical structures are obtained under investigation of these equations.

The properties of conservation laws are at the basis of the process of 
physical structure emergence. Therefore it is necessary to call 
attention to some properties and peculiarities of conservation laws.

\subsection{Properties and peculiarities of conservation laws.}

From the closure condition of exterior form it follows that the closed
inexact differential form is a conservative quantity
on some pseudostructure. From this one can see that the closed
inexact exterior differential form can correspond to conservation
law. The conservation laws for physical fields are just such
conservation laws.
{\footnotesize [The physical fields are a special form of the
substance, they are carriers of various interactions such as
electromagnetic, gravitational, wave, nuclear and other kinds of
interactions. The conservation laws for physical fields are
those that claim the existence of conservative physical quantities or
objects. Such conservation laws can be named the exact conservation
laws.]}

One can see that Inv. Structures made up by closed inexact form and 
relevant dual form correspond to conservation laws for physical fields.

The evolutionary skew-symmetric forms, from which, as it has been shown, the
closed inexact forms are obtained, correspond to conservation laws as 
well. However, these are conservation laws for material systems 
(material media). In contrast to conservation laws for physical fields, 
they are balance conservation laws (they establish the balance between 
the variation of physical quantities and external actions to the system) 
and are described by differential equations.
{\footnotesize [Material system is a variety of elements that have
internal structure and interact to one another. As examples of material
systems it may be thermodynamic, gas dynamical, cosmic systems, systems
of elementary particles and others]}.

The conservation laws for material systems are conservation laws
for energy, linear momentum, angular momentum, and mass.

The invariant structures corresponding to physical fields are obtained
from the equations that describe balance conservation laws for material 
media.

\subsection*{Analysis of the equations of conservation laws for material
systems.}

The balance conservation laws for energy, linear momentum, angular
momentum, and mass are described by partial differential equations [6].
(On examination of the integrability of these equations it is obtained 
the nonidentical relation that includes evolutionary form. From such 
evolutionary form the closed inexact forms and invariant structures 
corresponding to physical structures are obtained.)

Let us analyze the equations that describe the balance conservation
laws for energy and linear momentum.

In the accompanying frame of reference (this system is connected with 
the manifold made up by the trajectories of material system elements) 
the equations for energy and linear momentum are written as
$$
{{\partial \psi }\over {\partial \xi ^1}}\,=\,A_1 \eqno(16)
$$
$$
{{\partial \psi}\over {\partial \xi^{\nu }}}\,=\,A_{\nu },\quad \nu \,=\,2,\,...\eqno(17) 
$$

Here $\psi$  is the functional specifying the state of material system
(the action functional, entropy, wave function can be regarded as
examples of such a functional), $\xi^1$ is the coordinate along the
trajectory, $\xi ^{\nu }$ are the coordinates in the direction normal to
trajectory, $A_1$ is the quantity that depends on specific features
of material system and on external energy actions onto the system, and
$A_{\nu }$ are the quantities that depend on specific
features of material system and on external force actions.

Eqs. (16) and (17) can be convoluted into the relation
$$
d\psi\,=\,A_{\mu }\,d\xi ^{\mu },\quad (\mu\,=\,1,\,\nu )\eqno(18)
$$
where $d\psi $ is the differential
expression $d\psi\,=\,(\partial \psi /\partial \xi ^{\mu })d\xi ^{\mu }$. 

Relation (18) can be written as
$$
d\psi \,=\,\omega \eqno(19)
$$
here $\omega \,=\,A_{\mu }\,d\xi ^{\mu }$ is the skew-symmetric
differential form of first degree.

The relation obtained is an evolutionary relation.

Relation (19) was obtained from the equations of  balance conservation
laws for energy and linear momentum. In this relation the form $\omega $
is that of first degree. If the equations of balance conservation
laws for angular momentum be added to the equations for energy and 
linear momentum, this form in the evolutionary relation will be a form 
of second degree. And in  combination with the equation of balance
conservation law for mass this form will be the form of degree 3.

Thus, in general case the evolutionary relation can be written as
$$
d\psi \,=\,\omega^p \eqno(20)
$$
where the form degree  $p$ takes the values $p\,=\,0,1,2,3$..
(The evolutionary relation for $p\,=\,0$ is similar to that in
differential forms, and it was obtained from the interaction of energy
and time.)

The relations (19) and (20) are nonidentical evolutionary relations.

Let us show that the relation obtained from the equations
of balance conservation laws proves to be nonidentical.

To do so we shall analyze relation (19).

In the left-hand side of relation (19) there is the
differential that is a closed form. This form is an invariant
object. The right-hand side of relation (20) contains the differential
form $\omega$, which is not an invariant object since in real processes,
as it will be shown below, this form proves to be unclosed.
The commutator of this form is nonzero. The components of commutator of 
the form $\omega \,=\,A_{\mu }d\xi ^{\mu }$
can be written as follows:
$$
K_{\alpha \beta }\,=\,\left ({{\partial A_{\beta }}\over {\partial \xi ^{\alpha }}}\,-\,
{{\partial A_{\alpha }}\over {\partial \xi ^{\beta }}}\right )
$$
(here the term  connected with the manifold metric form
has not yet been taken into account).

The coefficients $A_{\mu }$ of the form $\omega $ have been obtained 
either from the equation of balance conservation law for energy or from 
that for linear momentum. This means that in the first case the 
coefficients depend on the energetic action and in the second case they 
depend on the force action. In actual processes energetic and force 
actions have different nature and appear to be inconsistent. 
The commutator of the form $\omega $ made up
of the derivatives of such coefficients is nonzero.
This means that the differential of the form $\omega $
is nonzero as well. Thus, the form $\omega$ proves to be unclosed and
cannot be a differential like the left-hand side.

This means that relation (19), as well as relation (20), cannot be 
identical ones. In such a way it can be shown that relation (20) is 
nonidentical as well.

Thus, the nonidentity of evolutionary relation
means that the balance conservation law equations are
inconsistent. And this indicates that the balance conservation laws are
noncommutative. (If the balance conservation laws be commutative,
the equations would be consistent and the evolutionary relation would
be identical).

The noncommutativity of balance conservation laws is a moving force of
evolutionary processes that proceed in material medium and lead to 
emergence of physical structures. This follows from the further analysis 
of the equations of balance conservation laws. The invariant structures 
obtained from these equations correspond to such physical structures.

\subsection{Mechanism of generation of physical structures.}

The relation obtained from the equations of balance
conservation laws involves the functional that specifies the
material system state. However, since this relation turns out to
be not identical, from this relation one cannot get the
differential $d\psi $  that could point out to the equilibrium state of
material system. The absence of differential means that the system
state is nonequilibrium. That is, in material system the internal
force acts.

As it has been already shown, the nonidentical evolutionary
relation turns out to be a selfvarying relation.

Selfvariation of the nonidentical evolutionary relation points to the
fact that the nonequilibrium state of material system turns out to be
selfvarying.
It is evident that this selfvariation proceeds under the action of
internal force whose quantity is described by commutator of the unclosed
evolutionary form $\omega^p $. (If the commutator
be zero, the evolutionary relation would be identical, and this would
point to the equilibrium state, i.e. the absence of internal forces.)
Everything that gives the contribution into the commutator of the form
$\omega^p $ leads to emergence of internal force.

Above it has been shown that under degenerate transformation from 
nonidentical evolutionary relation it can be obtained the identical 
relation
$$
d_\pi\psi=\omega_\pi^p\eqno(22)
$$

From such a relation one can obtain the
state function and this corresponds to equilibrium state of the system.
But identical relation can be realized only on pseudostructure (which
is specified by the condition of degenerate transformation). This
means that the transition of material system to equilibrium state
proceeds only locally. In other words, it is realized the transition
of material system from nonequilibrium state to locally
equilibrium one. In this case the total state of material system
remains to be nonequilibrium. The conditions of degenerate 
transformation can be caused by the degrees of freedom of material 
system.

As one can see from the analysis of nonidentical evolutionary
relation, the transition of material system from nonequilibrium
state to locally-equilibrium state proceeds spontaneously in the
process of selfvarying nonequilibrium state of material system
under realization of any degrees of freedom of this system.
(Translational degrees of freedom, internal degrees of freedom of
the system elements, and so on can be examples of such degrees of
freedom).

As it has been already said above, the transition from nonidentical
relation (21) obtained from balance conservation laws to identical
relation (22) means the following.
Firstly, the existence of state differential (left-hand side
of relation (22)) points to the transition of material system from
nonequilibrium state to locally-equilibrium state. And, secondly,
the emergence of closed (on pseudostructure) inexact exterior form
(right-hand side of relation (22)) points to the origination of physical
structure.
(Physical structures that are generated by material systems made up
physical fields.)

Thus one can see that the transition of material system from
nonequilibrium state to locally-equilibrium state is accompanied
by originating the differential-geometrical structures, which are
physical structures. The emergence of physical structures in the
evolutionary process reveals in material system as the emergence of
certain observable formations that develop spontaneously. Such
formations and their manifestations are fluctuations, turbulent
pulsations, waves, vortices, creating massless particles and others.
The intensity of such formations is controlled by a quantity
accumulated by the evolutionary form commutator at the instant in
time of originating physical structures.
The transition from evolutionary forms to closed exterior forms
describes such processes like the emergence of waves, vortices, turbulent
pulsations, the origination of massless particles and others [7].

Since the closed exterior forms corresponding to physical structures are
obtained from the evolutionary forms describing material systems, the
characteristics of physical structures are determined by characteristics of
material system generating these structures, and this enables one to 
classify physical structures by the parameters of
evolutionary forms and closed exterior forms.

As it has been shown above, the type of differential-geometrical 
invariant structures, and hence of physical structures (and,
accordingly, of physical fields) generated by the evolutionary
relation, depends on the degrees of differential form $p$ and $k$
and on the dimension of original inertial space $n$ (here $p$ is
the degree of evolutionary form in nonidentical relation
that is connected with the number of interacting balance
conservation laws, and $k$ is the  degree of closed form
generated by nonidentical relation). Introducing the
classification by numbers $p$, $k$, $n$ one can understand the
internal connection between various physical fields.

\bigskip

The above described mechanism of generation of physical structures
discloses an unique role of invariant structures in mathematical physics
and field theory.

It should be emphasized that such results were obtained due to using
the skew-symmetric exterior and evolutionary differential forms.
The mathematic apparatus of evolutionary forms, which describes the 
process of realization of closed exterior forms and invariant structures,
enables one to investigate the integrability of differential equations
(the conjugacy of the differential equations elements), discloses the 
mechanism of evolutionary processes, discrete transitions, quantum steps,
transitions from nonconjugated operators to conjugated ones,
and generation of various structures.
There are no such possibilities in any mathematical formalism.

1. Cartan E., Les Systemes Differentials Exterieus ef Leurs Application
Geometriques. -Paris, Hermann, 1945.

2. Bott R., Tu L.~W., Differential Forms in Algebraic Topology. 
Springer, NY, 1982.

3. Schutz B.~F., Geometrical Methods of Mathematical Physics. Cambridge
University Press, Cambridge, 1982.

4. Tonnelat M.-A., Les principles de la theorie electromagnetique 
et la relativite. Masson, Paris, 1959.

5. Petrova L.~I., Specific features of differential equations of 
mathematical physics. http://arxiv.org/abs/math-ph/0702019, 2007. 

6. Clark J.~F., Machesney ~M., The Dynamics of Real Gases. Butterworths, 
London, 1964. 

7. Petrova L.~I., The mechanism of generation of physical structures. 
//Nonlinear Acoustics - Fundamentals and Applications (18th International 
Symposium on Nonlinear Acoustics, Stockholm, Sweden, 2008) - New York, 
American Institute of Physics (AIP), 2008, pp.151-154.  

\end{document}